\documentclass[11pt, a4paper]{article}

\usepackage[utf8]{inputenc}
\usepackage[T1]{fontenc}
\usepackage{microtype, upref}
\usepackage[dvipsnames]{xcolor}

\usepackage{geometry}
\usepackage{amsmath, amssymb, amsthm, mathtools}
\usepackage{mathrsfs}
\usepackage{thmtools, thm-restate}

\usepackage{booktabs, longtable, tabularx}
\usepackage{adjustbox}
\usepackage{enumitem}

\usepackage{tikz}
\usetikzlibrary{arrows.meta,positioning,calc}

\usepackage[affil-it]{authblk}
\usepackage{setspace}
\allowdisplaybreaks

\usepackage[%
  colorlinks      = true,
  linktocpage     = true,
  linkcolor       = blue,
  citecolor       = blue,
  urlcolor        = blue,
  breaklinks      = true,
  hypertexnames   = false,
]{hyperref}
\usepackage{cleveref}

\usepackage[absolute,overlay]{textpos}

\newcommand{\toappear}[1]{%
  \AtBeginDocument{%
    \begin{textblock*}{150mm}(10mm,10mm) % (x,y) from page top-left
      {\small #1}%
    \end{textblock*}%
  }%
}

\newtheorem{theorem}{Theorem}[section]

\theoremstyle{definition}
\newtheorem{definition}[theorem]{Definition}
\newtheorem{remark}[theorem]{Remark}

\newcommand{\N}{\mathbb{N}}
\newcommand{\R}{\mathbb{R}}
\newcommand{\E}{\mathbb{E}}

\newcommand{\Var}{\mathrm{Var}}

\newcommand{\inner}[2]{\langle #1,#2\rangle}
\newcommand{\norm}[1]{\left\lVert #1\right\rVert}
\newcommand{\dd}{\,\mathrm{d}}

\newcommand{\GN}{G_N}
\newcommand{\Z}{\mathsf{Z}}
\newcommand{\FN}{\mathsf{F}_N}
\newcommand{\HN}{H_N}
\newcommand{\qstar}{q^{\ast}}

\newcommand{\sech}{\mathrm{sech}}
\newcommand{\inn}[1]{\langle #1 \rangle}

\title{\textbf{Michel Talagrand and the Rigorous Theory of Mean Field Spin Glasses}}
\author{Sourav Chatterjee\thanks{Department of Statistics, Stanford University, USA. Email: \href{mailto:souravc@stanford.edu}{\tt souravc@stanford.edu}.}}
\affil{Stanford University}
\date{} % (Abel Prize volume: omit running date)

\toappear{To appear in H.~Holden, R.~Piene (eds.): \emph{The Abel Prize 2023--2027}, Springer.}

\begin{document}
\maketitle

\begin{abstract}
Michel Talagrand played a decisive role in the transformation of mean-field spin glass theory into a rigorous mathematical subject.
This chapter offers a narrative account of that development.
We begin with the physical origins of the Sherrington--Kirkpatrick (SK) model and the emergence of the TAP and Almeida--Thouless stability frameworks, culminating in Parisi's replica symmetry breaking (RSB) ansatz and its hierarchical order parameter.
We then review early rigorous milestones, including high-temperature results and stability identities, and describe the consolidation of interpolation and cavity methods through the work of Guerra and of Aizenman--Sims--Starr.

The central event in this narrative is Talagrand's 2006 proof of the Parisi formula for the SK model and for a broad class of mixed \(p\)-spin models, and his subsequent analysis of Parisi measures.
We also discuss Talagrand's later program constructing pure states under extended Ghirlanda--Guerra identities and an atom at the maximal overlap, together with the structural results that followed, notably Panchenko's ultrametricity theorem and extensions of the Parisi formula.
Throughout, we indicate how related contributions by many authors fit into the same long-running program across probability, analysis, and mathematical physics.
\end{abstract}

%\tableofcontents

% ============================================================

% ============================================================
\section{Introduction}

This chapter is written for the Abel Prize volume honoring Michel Talagrand.
Its purpose is not merely to record a sequence of results, but to convey how Talagrand helped convert a set of extraordinarily detailed physical predictions into a durable mathematical theory:
one with canonical order parameters, stable proof architectures, and a shared language.
Accordingly, the emphasis here is on the evolution of ideas and methods, with Talagrand's contributions presented within the broader program they shaped and enabled.

Talagrand's entry into spin glasses came after foundational work in probability and analysis.
This included concentration and isoperimetric inequalities in product spaces (see, e.g., \cite{Talagrand1995Concentration,Talagrand1996NewConc}).
In the spin glass setting, those sensibilities---quantitative inequalities, stability under perturbations, and insistence on the correct order parameters---became organizing principles.
They are visible throughout the modern theory, from high-temperature control and interpolation bounds to the analysis of Parisi minimizers and the extraction of geometric structure from overlap identities.

The scope of mean field spin glass theory makes it impossible to be exhaustive in a single chapter.
Instead, the presentation is organized chronologically, with two recurring themes:
\begin{enumerate}
\item \textbf{From thermodynamics to geometry.}
The subject begins with the free energy, but it ultimately aims to describe the geometry of the Gibbs measure:
pure states, overlap distributions, and hierarchical organization (see \Cref{sec:mean-field-language,sec:identities-geometry,sec:pure-states}).
\item \textbf{From physical order parameters to probabilistic objects.}
Parisi's functional order parameter becomes, in the rigorous theory, a probability measure on \([0,1]\)
and, under suitable hypotheses, the law of overlaps among replicas (see \Cref{sec:parisi-functional,sec:parisi-measures}).
\end{enumerate}

\begin{remark}[On scope and attribution]
Spin glass theory is a large community, and no survey can cite every meaningful contribution.
The aim here is to be representative, emphasizing the main threads that intersect Talagrand's program:
interpolation and variational principles, stability identities, ultrametricity and cascades, universality, TAP and algorithmic connections, extremes/complexity, and related models.
\end{remark}

% ============================================================
\section{Physical origins: SK, TAP, AT, and Parisi  (1970s--1980s)}\label{sec:physics-prehistory}

This section sketches the physics prehistory that set the agenda for the rigorous theory.
We begin with the passage from disordered magnets and the Edwards--Anderson model to the mean field Sherrington--Kirkpatrick model (\Cref{sec:disordered-magnets}).
We then review the TAP approach and the replica-symmetric picture, and explain how stability considerations (notably the Almeida--Thouless criterion) led to the RS/RSB divide (\Cref{sec:TAP-AT}).
Finally, we outline Parisi's hierarchical ansatz and the conceptual shift from thermodynamics to overlap geometry that will reappear throughout the rest of the chapter.

\subsection{From disordered magnets to mean field}\label{sec:disordered-magnets}

A \emph{magnet} is a material whose microscopic magnetic moments (in simplified models, ``spins'') interact so that, at low enough temperature,
they may align and produce macroscopic magnetization.
A classical toy model is the Ising ferromagnet \cite{Ising1925}, where spins \(\sigma_i\in\{\pm1\}\) sit on a lattice and prefer to align:
\begin{equation*}
H^{\mathrm{Ising}}(\sigma)= -J\sum_{\langle i,j\rangle}\sigma_i\sigma_j - h\sum_i\sigma_i,
\end{equation*}
with \(J>0\) favoring agreement and \(h\) an external field, and $\langle i,j\rangle$ means that $i,j$ are neighbors on the lattice.

The Gibbs viewpoint is that, at inverse temperature \(\beta\), a configuration \(\sigma\) is observed at thermal equilibrium with probability proportional to \(e^{-\beta H(\sigma)}\),
so the normalizing constant
\(
\Z=\sum_{\sigma}e^{-\beta H(\sigma)}
\)
(the \emph{partition function}) encodes which macroscopic behaviors are typical.
A useful heuristic is to group configurations by the value of \(H(\sigma)\).
If \(L_a\) denotes the number of configurations with \(H(\sigma)\approx a\), then
\[
\mathrm{Prob}(H(\sigma)\approx a)\approx \Z^{-1}L_a e^{\beta J a} = \Z^{-1}e^{\beta J a + \log L_a}.
\]
For large systems this mass is typically sharply peaked near the maximizer of
\(
\beta J a+\log L_a,
\)
exhibiting the classic \emph{energy--entropy competition} and leading to (approximate) \emph{equivalence of ensembles}.
Differentiating \(\log \Z\) recovers typical values: 
\[
\partial_\beta \log \Z=-\langle H\rangle,
\]
where $\inn{\mathcal{O}}$ denotes the expected value of an observable $\mathcal{O}$, defined as
\[
  \inn{\mathcal{O}} = \frac{1}{\Z}\sum_\sigma \mathcal{O}(\sigma) e^{-\beta H(\sigma)}.
\]
A \emph{disordered magnet} is one in which the effective interaction strengths vary randomly in space.
Experimentally, this occurs for example in dilute magnetic alloys, where magnetic impurities are embedded in a non-magnetic host.
Such systems display a ``freezing'' transition without conventional long-range order, motivating the term \emph{spin glass};
see early experiments such as \cite{CannellaMydosh1972} and the classic experimental review \cite{Mydosh1993,BinderYoung1986}.

A minimal lattice model incorporating this disorder is the Edwards--Anderson (EA) model \cite{EdwardsAnderson1975}:
\begin{equation*}
H^{\mathrm{EA}}(\sigma)= -\sum_{\langle i,j\rangle} J_{ij}\sigma_i\sigma_j - h\sum_i\sigma_i,
\qquad \sigma\in\{\pm1\}^{\Lambda},
\end{equation*}
where \((J_{ij})\) are i.i.d.\ random couplings (often symmetric about \(0\)) and $\Lambda$ is a finite subset of $\mathbb{Z}^d$ for some $d$.
The competing signs of the \(J_{ij}\)'s lead to \emph{frustration}: it is impossible to satisfy all preferred pairwise alignments simultaneously.
This competition is the source of the complex low-temperature behavior that distinguishes spin glasses from ferromagnets.

Despite decades of work, the finite-dimensional EA model remains far less understood than its mean field counterpart.
There are multiple competing pictures for its low-temperature phase, and many basic structural questions remain open; see, e.g.,
the classic review \cite{BinderYoung1986} and the metastate/topical perspective \cite{NewmanStein2003},
as well as more recent broad reviews of spin-glass phenomenology and dynamics \cite{AltieriBaityJesi2023,DahlbergEtAl2025}.

The Sherrington--Kirkpatrick (SK) model \cite{SK1975} is the mean field analogue of EA:
every spin interacts with every other spin through a random coupling, scaled so that the total energy is extensive.\footnote{Here ``extensive'' means of order $N$: with the usual SK normalization, $H_N(\sigma)$ is typically of order  $N$ (so the energy per spin is $O(1)$), ensuring a nontrivial thermodynamic limit.}
In the pure SK (or pure 2-spin) model one takes i.i.d.\ standard Gaussian random variables \((g_{ij})_{1\le i<j\le N}\) and sets
\begin{equation}\label{eq:SK-informal}
H_N^{\mathrm{SK}}(\sigma)= -\frac{1}{\sqrt{N}}\sum_{1\le i<j\le N} g_{ij}\sigma_i\sigma_j
\qquad (\sigma\in\{-1,1\}^N),
\end{equation}
optionally adding an external field \(-h\sum_{i=1}^N\sigma_i\).
The mean field structure makes SK both more tractable and more ``universal'' as a mathematical object, but it also exposes new phenomena:
the Gibbs measure is predicted (and, in suitable senses, proved) to decompose into many competing ``phases'' (pure states),
with a nontrivial distribution of similarities (overlaps) between typical configurations.

The language that makes these statements precise appears later.
Overlaps and replica arrays are introduced in \Cref{sec:mean-field-language};
the ultrametric/hierarchical structure is discussed in \Cref{sec:identities-geometry};
and Talagrand's program of extracting pure states from overlap geometry is described in \Cref{sec:pure-states}.

\subsection{TAP, stability, and the RS/RSB divide}\label{sec:TAP-AT}

In the beginning, SK was analyzed via the \emph{replica method}: a formal approach to compute \(\E\log \Z_N\) (where $\Z_N$ is the partition function of the $N$-particle system) using the identity
\begin{equation}\label{eq:replica-trick}
\E\log \Z_N=\lim_{n\downarrow 0}\frac{\E \Z_N^n-1}{n},
\end{equation}
by first computing \(\E \Z_N^n\) for integers \(n\ge 1\) and then analytically continuing in \(n\) to \(0\). A common way physics presentations make the difficulty visible is to write the formal chain
\[
\E\log \Z_N
=\E\Big(\lim_{n\to 0}\frac{\Z_N^{\,n}-1}{n}\Big)=\lim_{n\to 0}\frac{\E \Z_N^{\,n}-1}{n} =\lim_{n\to 0}\frac{1}{n}\log \E \Z_N^{\,n},
\]
followed by an (even more delicate) exchange of \(\lim_{n\to 0}\) and \(\lim_{N\to\infty}\).
This is not a criticism of physics---it is a reminder of what the rigorous program had to replace by inequalities and limits.

Applied to SK, this produced the replica-symmetric (RS) ``solution'' of Sherrington and Kirkpatrick \cite{SK1975}.
However, already in physics it was clear that this solution becomes problematic at low temperature:
it can predict unphysical behavior (e.g., negative entropy in some parameter ranges), and it is unstable to certain fluctuations, as follows.

In the RS picture one posits that the macroscopic state is essentially described by a single number.
Concretely, the overlap between two independent replicas, \(R_{1,2}=N^{-1}\sum_i \sigma_i^1\sigma_i^2\), should concentrate near a deterministic value \(q\), often called the Edwards--Anderson order parameter. 
When RS fails, one expects that no single \(q\) can capture the equilibrium geometry, and the overlap distribution becomes nontrivial. De Almeida and Thouless \cite{AT1978} performed a stability analysis of the RS ansatz,
finding a boundary (the ``AT line'' in the \((\beta,h)\) plane) across which RS becomes unstable to replica symmetry breaking (RSB).
We discuss the RSB picture in the next subsection (\Cref{sec:parisi-physics}).
In the rigorous story, this stability condition reappears as a criterion for when the \emph{Parisi minimizer} is replica-symmetric
(see \Cref{sec:parisi-measures} for details).

A more mechanistic  mean field approach\footnote{Here ``mechanistic'' means that the approach aims to describe equilibrium by (approximately) solving self-consistency equations for local observables (such as the magnetizations $m_i=\langle \sigma_i\rangle$), rather than by directly postulating a global variational formula for the free energy.} was given by the Thouless--Anderson--Palmer (TAP) equations \cite{TAP1977},
which aimed to describe the Gibbs measure through its local magnetizations \(m_i=\langle\sigma_i\rangle\).
TAP yields a self-consistency relation of the form
\begin{equation}\label{eq:TAP-informal}
m_i \approx \tanh\!\Big(h + \sum_{j\ne i}\frac{\beta}{\sqrt{N}}g_{ij}m_j - \beta^2(1-q)\,m_i\Big),
\end{equation}
where the last term is the Onsager ``reaction'' correction and \(q\) is the typical self-overlap \(N^{-1}\sum_i m_i^2\) (which, in the replica-symmetric regime, turns out to match the typical value of $R_{1,2}$ discussed above).
In the high-temperature regime, TAP and RS match; beyond that regime, one expects many TAP solutions (reflecting many states).

\subsection{Parisi's replica symmetry breaking and a measure-valued order parameter}\label{sec:parisi-physics}

Parisi's 1979 breakthrough \cite{Parisi1979} can be read as a response to two pressures:
the RS replica computation is unstable (AT), and TAP suggests a profusion of metastable states.
Parisi proposed that one should \emph{break} replica symmetry hierarchically, producing a variational principle over order parameters.
In the physics literature this culminates in the Parisi free energy formula and a detailed geometric picture of the Gibbs measure.

Concretely, Parisi's idea was to replace the single RS overlap parameter \(q\) by a \emph{step function} \(q(x)\), \(x\in[0,1]\), encoding a nested (hierarchical) organization of replicas:
one partitions the replicas into blocks of size \(m_1\), then partitions each block into sub-blocks of size \(m_2\), and so on.
Replicas in the same smallest block are assigned the largest overlap value \(q_k\); replicas that separate earlier in the hierarchy are assigned smaller overlaps \(q_{k-1},q_{k-2},\dots\).
This produces an explicit ``\(k\)-step RSB ansatz'' parameterized by numbers
\[
0\le q_0\le q_1\le\cdots\le q_k\le 1,
\qquad
0=m_0<m_1<\cdots<m_k< m_{k+1}=1,
\]
and a corresponding formula for the free energy depending on these parameters.
Parisi then takes \(k\to\infty\), leading (in modern language) to a variational problem over probability measures on \([0,1]\) (equivalently, distribution functions \(x\mapsto \alpha([0,x])\)), which is the rigorous \emph{Parisi functional} introduced later in \Cref{sec:parisi-functional}.
The synthesis by M\'ezard--Parisi--Virasoro \cite{MPV1987} developed three interlocking claims:
\begin{itemize}
\item the Gibbs measure decomposes into many pure states (clusters with internal coherence);
\item overlaps between replicas have a nontrivial distribution;
\item the joint overlap structure is \emph{ultrametric} (hierarchical).
\end{itemize}
We state a rigorous form of the Parisi free energy variational principle in \Cref{sec:parisi-functional}
(and a discussion of Talagrand's proof in \Cref{sec:parisi-functional}).
The overlap language needed to formulate the ``geometry'' is introduced in \Cref{sec:mean-field-language},
and the ultrametric/cascade picture is explained in \Cref{sec:identities-geometry}.

\subsection{Solvable hierarchies: REM/GREM and Ruelle cascades}\label{sec:rem-grem}

The Parisi picture is conceptually powerful but, at the level of SK itself, it was for a long time supported mainly by physics arguments.
A complementary development was therefore to study models where the same hierarchical ideas can be proved directly, and where the limiting Gibbs weights and overlap structure can be computed explicitly.

Derrida introduced the random energy model (REM) \cite{Derrida1981} as a solvable caricature of the freezing transition:
energies of different configurations are i.i.d., so the partition function becomes an extreme-value object.
The generalized REM (GREM) correlates energies through a tree, making hierarchy explicit and producing nontrivial overlap distributions.
These models served as a laboratory for the emerging RSB ideas: one can compute the free energy, see ``many states,''
and verify that the relevant weights are Poisson--Dirichlet in the glassy phase.

Ruelle gave a mathematical formulation of these hierarchies via probability cascades \cite{Ruelle1987},
which later became canonical limiting objects for ultrametric Gibbs measures.
Bolthausen--Sznitman \cite{BolthausenSznitman1998} further developed the cascade viewpoint in probabilistic terms.

These solvable hierarchical models matter for SK for two reasons. First, they make Parisi's hierarchical ansatz concrete and provable.
Second, they provide the probabilistic templates (Poisson--Dirichlet weights, cascades) that later reappear
in Talagrand's pure-state program (\Cref{sec:pure-states}) and in structural theorems driven by overlap identities (\Cref{sec:identities-geometry}).

\subsection{A brief roadmap of the physics predictions}\label{sec:physconj}
To keep the narrative anchored, it is helpful to record the core claims of the physics theory as concrete mathematical statements.
They will reappear later as theorems (or as theorems under hypotheses). For the SK model (and more generally mixed \(p\)-spin models), physics predicted:
\begin{enumerate}
\item[(i)] \textbf{Parisi variational principle (free energy).}
The limit \(\lim_{N\to\infty}\FN(\beta,h)\) exists and equals a variational value \(\mathcal{P}(\xi,h)\) over probability measures on \([0,1]\)
(defined precisely in \Cref{sec:parisi-functional}).
\item[(ii)] \textbf{Replica symmetry vs.\ symmetry breaking.}
In the RS region, the overlap \(R_{1,2}\) concentrates near a deterministic value \(q\), and the free energy equals the RS functional at \(q\).
In the RSB region, the overlap distribution is nontrivial, with many overlap levels.
\item[(iii)] \textbf{Ultrametricity (geometry of overlaps).}
In the RSB region, the joint overlap array of many replicas is ultrametric, in the sense that for typical triples
\(R_{1,2}\ge \min\{R_{1,3},R_{2,3}\}\).
\item[(iv)] \textbf{Pure states and cascades.}
The Gibbs measure decomposes into (infinitely many) pure states with random weights, and these weights exhibit Poisson--Dirichlet/cascade statistics.
\end{enumerate}
Historically, these were derived non-rigorously.
In the rigorous story, (i) becomes Talagrand's 2006 theorem (\Cref{sec:parisi-functional}) and its extensions; (ii) becomes the characterization of the replica-symmetric region in terms of the Parisi minimizer (a one-atom Parisi measure, equivalently an RS minimizer at some $q$, yielding concentration of $R_{1,2}$) and its breakdown in the RSB region (a nontrivial Parisi measure and hence a nontrivial overlap distribution with multiple overlap levels; see \Cref{sec:parisi-measures} for the Parisi-measure viewpoint and \Cref{sec:TAP-AT} for the physics RS/RSB divide); (iii) becomes Panchenko's ultrametricity theorem under overlap identities (\Cref{sec:identities-geometry}); and (iv) is connected to Talagrand's 2008--2010 pure-state program (\Cref{sec:pure-states}).

% ============================================================
\section{Mean field models and overlaps: the common language}\label{sec:mean-field-language}

This section fixes notation and collects the basic objects that appear throughout the rigorous theory.
We define the Gibbs measure, free energy, and replicas, and then describe the mixed \(p\)-spin family via its covariance function \(\xi\).
Later sections use these definitions to formulate stability identities, interpolation bounds, and the Parisi variational principle.

\subsection{Gibbs measures, replicas, and overlaps}

Throughout, \(\Sigma_N=\{-1,1\}^N\). 
For the SK model introduced in \Cref{sec:disordered-magnets}, we define the Hamiltonian slightly differently than earlier, by absorbing the inverse temperature  into the Hamiltonian itself: 
\begin{equation}\label{eq:SK}
\HN(\sigma)= -\frac{\beta}{\sqrt{N}}\sum_{1\le i<j\le N} g_{ij}\sigma_i\sigma_j,
\qquad \sigma\in\Sigma_N.
\end{equation}

\paragraph{Gibbs measure and free energy.}
In addition to the Hamiltonian defined above, we insert an external magnetic field term $- h\sum_{i=1}^N \sigma_i$. After this insertion, 
the partition function and Gibbs measure are
\begin{equation}
\label{eq:ZN-GN}
\begin{split}
&\Z_N(\beta,h)=\sum_{\sigma\in\Sigma_N}\exp\Bigl(-\HN(\sigma)- h\sum_{i=1}^N \sigma_i\Bigr), \\ 
&\GN(\sigma)=\frac{1}{\Z_N(\beta,h)}\exp\Bigl(-\HN(\sigma)- h\sum_{i=1}^N \sigma_i\Bigr).
\end{split}
\end{equation}
The expected free energy per spin is
\begin{equation}\label{eq:FN}
\FN(\beta,h)=\frac{1}{N}\E\log \Z_N(\beta,h).
\end{equation}

\paragraph{Replicas and their relation to the replica trick.}
In the physics literature, \emph{replicas} first appear in the replica trick \eqref{eq:replica-trick}, where one formally studies \(\E \Z_N^n\) for integer \(n\) by introducing \(n\) copies of the spin configuration and then (non-rigorously) extrapolating to \(n\downarrow 0\).
In the rigorous theory, the same word is used in a more literal way: a \emph{replica} is simply an independent sample from the Gibbs measure \(\GN\).
The bridge between these viewpoints is that computations of \(\E \Z_N^n\) involve \(n\) spins sampled from the same disorder, and their \emph{mutual overlaps} are exactly the quantities that survive (and organize) the calculation.

\begin{definition}[Overlap]\label{def:overlap}
For independent replicas \(\sigma^1,\sigma^2\sim \GN\),
\begin{equation}\label{eq:overlap}
R_{1,2}=\frac{1}{N}\sum_{i=1}^N \sigma_i^1\sigma_i^2\in[-1,1].
\end{equation}
For \(s\) replicas, define the overlap array \(\big(R_{\ell,\ell'}\big)_{1\le \ell<\ell'\le s}\).
\end{definition}

\paragraph{Why overlaps are the natural order parameter.}
For mean field spin glasses the disorder is Gaussian, so the model is (in law) determined by covariances of the Hamiltonian.
In SK with zero external field, a short computation shows that for two configurations \(\sigma^1,\sigma^2\),
\begin{equation}\label{eq:SK-cov}
\E\big[\HN(\sigma^1)\HN(\sigma^2)\big]
=\frac{\beta^2}{N}\sum_{1\le i<j\le N}(\sigma_i^1\sigma_i^2)(\sigma_j^1\sigma_j^2)
=\frac{\beta^2 N}{2}\Big(R_{1,2}^2-\frac{1}{N}\Big),
\end{equation}
so the covariance depends on \((\sigma^1,\sigma^2)\) only through their overlap \(R_{1,2}\) (up to a negligible \(O(1)\) correction).
This is the basic reason overlap arrays are the right state variables in interpolation arguments and in the structural theory.

\subsection{Mixed \(p\)-spin models and the covariance function \(\xi\)}\label{sec:mixed-p-spin}

Historically, \(p\)-spin interaction models were introduced in the physics literature as mean field models with genuine many-body interactions,
and (especially for \(p\ge 3\)) as tractable toy models for structural glasses and complex energy landscapes;
see, e.g., \cite{GrossMezard1984,Gardner1985,KirkpatrickThirumalai1987}.
From the rigorous viewpoint, allowing \emph{mixtures} of different \(p\)'s is not just generality for its own sake:
mixtures are stable under small perturbations, they allow one to enforce generic overlap identities,
and many fundamental arguments (Guerra's interpolation, Talagrand's proof strategy) are naturally expressed in terms of a single covariance function \(\xi\).

\begin{definition}[Mean field Gaussian model via \(\xi\)]\label{def:xi-model}
A broad class of models has centered Gaussian Hamiltonians with covariance
\begin{equation}\label{eq:xi-cov}
\E\big[\HN(\sigma^1)\HN(\sigma^2)\big]=N\,\xi(R_{1,2}),
\end{equation}
where \(\xi:[-1,1]\to\R\) is even and typically \(\xi(q)=\sum_{p\ge 1}\beta_p^2 q^p\) with \(\sum_p \beta_p^2<\infty\).
\end{definition}
A standard realization of \eqref{eq:xi-cov} is the mixed \(p\)-spin Hamiltonian
\begin{equation}\label{eq:mixed-p-spin}
\HN(\sigma)
=
-\sum_{p\ge 1}\frac{\beta_p}{N^{(p-1)/2}}
\sum_{1\le i_1,\dots,i_p\le N} g_{i_1,\dots,i_p}\,\sigma_{i_1}\cdots\sigma_{i_p},
\end{equation}
where the \(g_{i_1,\dots,i_p}\) are i.i.d.\ standard Gaussians (often symmetrized in the indices).
A short computation yields \(\E[\HN(\sigma^1)\HN(\sigma^2)]=N\sum_{p\ge1}\beta_p^2R_{1,2}^p\), i.e., \eqref{eq:xi-cov}.

The SK model corresponds to \(\xi(q)=\tfrac{\beta^2}{2}q^2\).
A major reason this formulation matters is that many theorems (including Guerra's bound and Talagrand's proof strategy)
are naturally stated in terms of \(\xi\), and convexity/regularity properties of \(\xi\) control monotonicity
in interpolation arguments.

% ============================================================
\section{Early rigorous anchors (1980s--1990s): high temperature and stability identities}\label{sec:early-anchors}

Before the transformative breakthroughs by Talagrand, several important rigorous pillars were already in place.
Seen through the lens of the physics conjectures discussed in \Cref{sec:physconj}, these results provided the first firm evidence for the two central lessons:
(i) there is a genuinely ``classical'' (RS) high-temperature regime, and
(ii) overlap observables in the low-temperature regime cannot behave like ordinary self-averaging\footnote{Here ``self-averaging'' means concentration with respect to the disorder: a sequence of random observables $X_N$ is self-averaging if its disorder-to-disorder fluctuations vanish asymptotically, e.g., $\Var(X_N)\to 0$ (equivalently, $X_N$ converges in probability to a deterministic limit).}   quantities. 
They also introduced the algebraic overlap identities that later became the engine behind ultrametricity and pure-state structure.

In this section we sketch this pre-Talagrand rigorous landscape: first, high-temperature results where overlaps behave in a replica-symmetric way and the free energy can be computed sharply; then, low-temperature indicators that a single deterministic order parameter cannot self-average; and finally, the stability/consistency identities (Aizenman--Contucci and Ghirlanda--Guerra) that later drive the geometric theory.

\paragraph{High temperature: the free energy is essentially annealed (Aizenman--Lebowitz--Ruelle).}
Aizenman, Lebowitz, and Ruelle \cite{AizenmanLebowitzRuelle1987} proved foundational results for SK in high temperature,
and Comets and Neveu \cite{CometsNeveu1995} developed related analyses using stochastic calculus.

\begin{theorem}[High-temperature SK at \(h=0\) (representative form; cf.\ \cite{AizenmanLebowitzRuelle1987,CometsNeveu1995})]\label{thm:ALR-highT}
Assume \(h=0\) and \(\beta<1\). Then
\[
\lim_{N\to\infty}\FN(\beta,0)=\log 2+\frac{\beta^2}{4},
\]
and overlaps are replica-symmetric in the sense that \(\E\langle R_{1,2}^2\rangle\to 0\).
Moreover, the fluctuations of \(\log \Z_N(\beta,0)\) around its mean are asymptotically Gaussian (in a normalization made precise in \cite{AizenmanLebowitzRuelle1987}).
\end{theorem}

\paragraph{Low temperature: the overlap order parameter cannot self-average (Pastur--Shcherbina).}
Pastur and Shcherbina \cite{PasturShcherbina1991} proved an influential theorem that ties correctness of the RS free energy to self-averaging of a natural order parameter.
To state it, let
\[
\bar q_N=\frac{1}{N}\sum_{i=1}^N \langle \sigma_i\rangle^2,
\]
where \(\langle\cdot\rangle\) is Gibbs expectation for fixed disorder.

\begin{theorem}[Pastur--Shcherbina: self-averaging would force the SK (RS) formula \cite{PasturShcherbina1991}]\label{thm:PS}
If \(\Var(\bar q_N)\to 0\) as \(N\to\infty\) (i.e., \(\bar q_N\) is self-averaging), then the limiting free energy would coincide with the Sherrington--Kirkpatrick replica-symmetric expression.
\end{theorem}

In modern hindsight, once the Parisi formula is known, one can read this theorem contrapositively:
in parameter regions where the Parisi value differs from the RS value (for example, at \(h=0\) and sufficiently large \(\beta\)),
the order parameter \(\bar q_N\) cannot be self-averaging.
This is an early rigorous signal that the low-temperature Gibbs measure cannot be described by a single macroscopic state.

\paragraph{Overlap identities as a structural engine (Aizenman--Contucci; Ghirlanda--Guerra).}
Two 1998 developments became central to the later structural theory.

\begin{itemize}
\item \textbf{Stochastic stability (Aizenman--Contucci) \cite{AizenmanContucci1998}.}
One perturbs the Hamiltonian by adding a small independent Gaussian field \(t \widetilde H_N\) (with the same covariance structure) and asks that, in the thermodynamic limit, the overlap statistics are stable under this perturbation.
Differentiating suitable expectations with respect to \(t\) at \(t=0\) produces linear identities among overlap moments (the Aizenman--Contucci identities).
\item \textbf{Ghirlanda--Guerra identities \cite{GhirlandaGuerra1998}.}
These relate mixed overlap moments in a way that enforces a remarkable self-consistency of the overlap array.
\end{itemize}

To highlight what the Ghirlanda--Guerra identities say, one common (simplified) form is:
for a bounded function \(f\) of the overlaps among \(s\) replicas and for any integer \(p\ge 1\),
\begin{equation}\label{eq:GG-schematic}
\E\big\langle f\,R_{1,s+1}^p\big\rangle
\approx
\frac1s\,\E\langle f\rangle\,\E\langle R_{1,2}^p\rangle
+
\frac1s\sum_{\ell=2}^s \E\big\langle f\,R_{1,\ell}^p\big\rangle,
\end{equation}
where \(\langle\cdot\rangle\) denotes the Gibbs expectation (and ``\(\approx\)'' hides an error that vanishes as \(N\to\infty\),
often after adding a tiny perturbation to enforce the identities).
The right-hand side says: a new replica overlaps with replica \(1\) either like an independent draw (first term) or by copying one of the existing overlaps (second term).
This ``copy-or-refresh'' structure is a key engine behind ultrametricity and cascades (\Cref{sec:identities-geometry})
and is one of the conceptual bridges from thermodynamics to geometry.

\begin{remark}
A recurring methodological move (used systematically later) is to add a tiny generic perturbation to the Hamiltonian.
This typically does not change the free energy limit, but it enforces strong overlap identities.
\end{remark}

% ============================================================
\section{Talagrand's early spin glass period and the field's reorganization (1998--2003)}\label{sec:talagrand-early}

Talagrand's entry into the SK problem around 1998 helped convert the subject from a ``physics solution in search of a proof'' into a sustained probability program.
This period also coincides with a reorganization of the field around a small set of decisive ideas:
overlaps as the right coordinates, cavity comparisons as a way to relate system sizes, and interpolation inequalities as the main engine for sharp bounds.
In this section we sketch how Talagrand's papers from 1998--2003 contributed to that reorganization and prepared the ground for the 2006 Parisi formula.

\subsection{From the challenge paper to quantitative methods (1998--2000)}\label{sec:talagrand-1998-2000}

We start with the 1998 ``challenge'' paper, which distilled the SK problem into a list of concrete mathematical tasks and insisted that the overlap array---not just the limiting free energy---is the right object to organize the theory.
We then turn to the 2000 papers that supplied quantitative tools (notably exponential tail bounds) and that stressed the flexibility of the Parisi picture by exhibiting models with arbitrarily deep finite-step replica symmetry breaking, before briefly noting low-temperature results for mean field $p$-spin models as a laboratory for RSB phenomena.

\paragraph{1998: what needed to be proved, and why overlaps are unavoidable.}
In ``\emph{The Sherrington--Kirkpatrick model: a challenge for mathematicians}'' \cite{Talagrand1998SKChallenge}, Talagrand framed SK for a mathematical audience in a way that made the obstructions explicit.
Two themes from that paper became lasting:
\begin{itemize}
\item The object is not merely \(\log \Z_N\), but the random Gibbs measure and the overlap statistics it induces.
\item One must replace the replica trick by robust inequalities, with cavity arguments (induction on \(N\)) as a natural organizing device.
\end{itemize}
The paper also contains concrete results and mechanisms.
For example, it proves a small-\(\beta\) replica-symmetric formula in nonzero external field (a precursor of \Cref{thm:RS-small-beta} below),
and it develops implications of (hypothetical) self-averaging for order parameters: roughly, ``if an order parameter self-averages, then the SK (RS) formula must hold,''
thereby turning failure of the RS formula into a rigorous route to non-self-averaging phenomena (compare \Cref{thm:PS}).

\paragraph{2000: exponential inequalities and models with many symmetry-breaking levels.}
Talagrand's 2000 papers strengthened the quantitative infrastructure and demonstrated the flexibility of the RSB picture.
In \cite{Talagrand2000RSBExp}, he developed sharp exponential inequalities for SK-type models; these are precisely the kinds of uniform controls
needed to make long cavity/interpolation arguments rigorous without losing track of error terms.
In \cite{Talagrand2000MultiLevels}, he showed that finite-step replica symmetry breaking is not a peculiarity of a single Hamiltonian.

To phrase this in accessible terms, recall (from \Cref{sec:physconj}) that RSB corresponds to multiple overlap levels.
A \emph{\(k\)-step RSB} order parameter means, informally, that the predicted overlap distribution is supported on \(k+1\) distinct values (equivalently: the cumulative order parameter has \(k\) genuine jumps).

\begin{theorem}[Models with arbitrarily deep finite-step RSB (informal; \cite{Talagrand2000MultiLevels})]\label{thm:finite-rsb-depth}
For each \(k\in\N\), there exist mixed \(p\)-spin models (equivalently: choices of the covariance function \(\xi\))
for which the Parisi minimizer has at least \(k+1\) support points.
Equivalently, one can realize models with at least \(k\) nontrivial symmetry-breaking steps.
\end{theorem}

Thus Talagrand's theorem shows that finite-step RSB depth can be made arbitrarily large,
but it does not, in general, provide an exact prescription of a preassigned finite number of steps.

In a third important 2000 paper, Talagrand developed low-temperature results for mean field \(p\)-spin models \cite{Talagrand2000pSpinLowTemp},
which serve as laboratories for RSB phenomena and influenced later structural intuition.

\subsection{Replica symmetry as theorem, and a book that codified the methodological framework (2001--2003)}\label{sec:talagrand-2001-2003}

The next step in this story is the emergence of replica symmetry as a theorem in the high-temperature regime: overlaps concentrate at a deterministic fixed point, and the free energy is given by the one-parameter RS variational formula.
We then briefly note how Talagrand's 2003 monograph helped codify the cavity and interpolation toolkit, giving the subject a standard proof architecture that it would rely on in the years leading up to the Parisi formula.

\paragraph{The replica-symmetric (high-temperature) regime.}
Talagrand's 2002 work \cite{Talagrand2002HighTemp} gives a careful treatment of the high-temperature (RS) phase of SK.
A useful way to phrase RS in rigorous terms is:
\begin{quote}
Replica symmetry means that the overlap \(R_{1,2}\) concentrates near a deterministic value \(q\),
and the free energy equals the replica-symmetric variational expression evaluated at \(q\).
\end{quote}
For SK, the RS candidate is described by the self-consistency equation
\begin{equation}\label{eq:RS-fixed-point}
q=\E\tanh^2(\beta\sqrt{q}\,Z+h),\qquad Z\sim\mathcal{N}(0,1),
\end{equation}
and the corresponding RS free energy functional is
\begin{equation}\label{eq:RS-free-energy}
\mathcal{P}_{\mathrm{RS}}(\beta,h;q)
=
\log 2 + \E\log\cosh(\beta\sqrt{q}\,Z+h) + \frac{\beta^2}{4}(1-q)^2.
\end{equation}
(Equivalently: \(\mathcal{P}_{\mathrm{RS}}\) is the Parisi functional restricted to one-atom order parameters.)

\begin{theorem}[Replica-symmetric formula for small \(\beta\) with external field (Talagrand~\cite{Talagrand2002HighTemp}, representative form)]\label{thm:RS-small-beta}
There exists \(\beta_0>0\) such that for every \(h>0\) and every \(0\le \beta\le \beta_0\),
the fixed point equation \eqref{eq:RS-fixed-point} has a unique solution \(q\in[0,1]\) and
\[
\lim_{N\to\infty}\FN(\beta,h)= \log 2 + \E\log\cosh(\beta\sqrt{q}\,Z+h) + \frac{\beta^2}{4}(1-q)^2,
\qquad Z\sim\mathcal{N}(0,1).
\]
Moreover, overlaps concentrate around \(q\).
\end{theorem}

\paragraph{Sketch of the method: cavity plus interpolation.}
The proof strategy that emerges clearly in Talagrand's work is a blend of two ideas.

(i) \emph{Cavity (induction on \(N\)).}
Compare the \(N\)-spin system to the \((N-1)\)-spin system by isolating the last spin \(\sigma_N\).
Conditioned on the disorder among the first \(N-1\) spins, the last spin sees an \emph{effective random field}
whose distribution is approximately Gaussian, with variance controlled by overlap-type quantities of the \((N-1)\)-spin Gibbs measure.
In high temperature, one can show these overlap quantities are tightly concentrated, so the effective field is close to a deterministic Gaussian mixture governed by a single parameter \(q\).

(ii) \emph{Interpolation/IBP (integration by parts; Gaussian calculus).}
Introduce an interpolation between the true \(N\)-spin Hamiltonian and a tractable reference system whose free energy increment can be computed explicitly in terms of \(q\).
Gaussian integration by parts identifies the derivative of the interpolated free energy as an average of overlap fluctuations.
In a high-temperature regime these fluctuations can be shown to be small (or the derivative can be shown to have the correct sign),
which closes the argument and forces the overlap to concentrate at the self-consistent solution of \eqref{eq:RS-fixed-point}.
Evaluating the free energy increment then yields \eqref{eq:RS-free-energy}.

This ``cavity + interpolation'' posture becomes one of the canonical templates of the rigorous theory.

\paragraph{2003: a synthesis and a pivot.}
Talagrand's 2003 monograph \cite{TalagrandBook2003} is a bridge book: it is written for mathematicians entering the subject,
it supplies complete proofs of much of what was then known, and it codifies what became standard technology.
Beyond the SK model itself, the book treats (among other topics) the Almeida--Thouless line and second-moment ideas,
Guerra's broken-replica-symmetry bound and the Parisi functional, the perceptron (in several settings),
the Hopfield model, and the \(p\)-spin interaction model at low temperature.
It also contains a striking feature that later readers came to associate with Talagrand's style:
a carefully curated set of open problems of varying difficulty, which helped define a research agenda while the rigorous theory was still rapidly evolving.

In two short 2003 notes \cite{Talagrand2003GuerraBound,Talagrand2003MeaningOrderParameter},
Talagrand sharpened the community's focus:
Guerra's new interpolation bound was the right inequality in the right language,
and the Parisi order parameter should be treated as the central object to interpret.
We explain Guerra's interpolation scheme and its resulting Parisi upper bound in \Cref{sec:guerra-bound}.

Around the same period, Talagrand also worked on disordered mean field models beyond SK,
including random \(k\)-SAT and the Hopfield model \cite{Talagrand2001KSAT,Talagrand2001HopfieldCritical}.
These papers helped establish that spin glass methods are  broadly applicable across disordered systems.

% ============================================================
\section{Interpolation and variational principles (2002–2003)}\label{sec:interpolation-revolution}

The period 2002--2003 reorganized SK around interpolation/cavity principles.
Guerra--Toninelli established the thermodynamic limit.
Guerra then produced the Parisi upper bound by a multilevel RSB interpolation.
Aizenman--Sims--Starr (ASS) recast the problem as a variational principle over overlap structures,
clarifying why Ruelle cascades are natural candidates.
These ideas provide the scaffolding on which Talagrand's 2006 proof rests.

\subsection{Thermodynamic limit}\label{sec:thermodynamic-limit}

Guerra and Toninelli established existence of the thermodynamic limit for a broad class of mean field spin glasses \cite{GuerraToninelli2002}.
In our setting (Gaussian mean field models with covariance \(N\xi(R_{1,2})\)), a representative formulation is:

\begin{theorem}[Thermodynamic limit (representative form) \cite{GuerraToninelli2002}]\label{thm:thermo-limit}
For mean field Gaussian spin glasses with covariance \(N\xi(R_{1,2})\) (with suitable regularity of \(\xi\)),
then
\(
\lim_{N\to\infty}\FN(\beta,h)
\)
exists.
Moreover, by standard Gaussian concentration (see, e.g., \cite{Talagrand1995Concentration,Talagrand1996NewConc} or the treatments in \cite{TalagrandBook2003,TalagrandMF1}),
one also has almost sure convergence of \(\frac1N\log\Z_N\) to the same limit.
\end{theorem}

This ensures there is a well-defined number for Parisi's functional to match.

\subsection{Guerra's RSB bound}\label{sec:guerra-bound}

Guerra introduced the replica symmetry breaking interpolation and proved the Parisi \emph{upper bound} \cite{Guerra2003BrokenRSB}.
At a high level, the scheme is:
\begin{itemize}
\item choose a discrete RSB ``grid'' \((m,q)\) (formalized in \Cref{sec:parisi-functional});
\item build a hierarchical Gaussian reference field whose covariances encode this grid;
\item interpolate between the true Hamiltonian and the hierarchical reference system using Gaussian integration by parts;
\item show that the derivative of the interpolated free energy has a sign controlled by convexity/positivity,
yielding a one-sided inequality whose endpoint is the Parisi functional evaluated at \((m,q)\).
\end{itemize}
Optimizing over all discrete grids then yields the Parisi bound in its measure-valued form.

\paragraph{The discrete Parisi functional (needed for the bound).}
We fix a standard discrete formulation (the one naturally produced by Guerra's multilevel interpolation).

\begin{definition}[Discrete Parisi parameters and functional]\label{def:discrete-parisi}
Fix \(k\in\N\) and sequences
\[
0=m_0\le m_1\le \cdots \le m_k=1,\qquad
0=q_0\le q_1\le \cdots \le q_{k+1}=1.
\]
Let \((z_p)_{p=0}^k\) be independent Gaussian random variables with
\[
\E z_p^2=\xi'(q_{p+1})-\xi'(q_p).
\]
Define
\[
X_{k+1}=\log\cosh\Big(h+\sum_{p=0}^k z_p\Big),
\]
and recursively
\[
X_\ell=
\begin{cases}
\frac{1}{m_\ell}\log \E_\ell \exp(m_\ell X_{\ell+1}), & m_\ell>0,\\[0.4em]
\E_\ell[X_{\ell+1}], & m_\ell=0,
\end{cases}
\]
where \(\E_\ell\) denotes expectation in \((z_p)_{p\ge \ell}\).
Let \(\theta(q)=q\xi'(q)-\xi(q)\) and define
\[
\mathcal{P}_k(m,q)=\log 2 + X_0 - \frac12\sum_{\ell=1}^k m_\ell\big(\theta(q_{\ell+1})-\theta(q_\ell)\big).
\]
Finally set \(\mathcal{P}(\xi,h)=\inf_{k,m,q}\mathcal{P}_k(m,q)\).
\end{definition}

\begin{theorem}[Guerra's RSB bound (Parisi upper bound) \cite{Guerra2003BrokenRSB}]\label{thm:guerra-bound}
For the SK model and, more generally, mean field Gaussian spin glasses with covariance \(N\xi(R_{1,2})\) (with \(\xi\) even and convex on \([0,1]\)),
for every choice of discrete Parisi parameters \((k,m,q)\) as in \Cref{def:discrete-parisi},
\[
\limsup_{N\to\infty}\FN(\beta,h)\le \mathcal{P}_k(m,q).
\]
Consequently,
\[
\limsup_{N\to\infty}\FN(\beta,h)\le \mathcal{P}(\xi,h).
\]
\end{theorem}

\paragraph{A proof sketch in one paragraph.}
Let \(\phi_N(t)=\frac{1}{N}\E\log \sum_{\sigma}\exp(-H_t(\sigma)+h\sum_i\sigma_i)\),
where \(H_t\) interpolates between the original Hamiltonian and a carefully chosen hierarchical Gaussian field
associated to the grid \((m,q)\).
By Gaussian integration by parts, \(\phi_N'(t)\) can be written explicitly in terms of overlap fluctuations of multiple replicas.
The hierarchical reference field is engineered so that these fluctuations enter through a nonnegative quadratic form (or a form with the correct sign),
which implies \(\phi_N'(t)\le 0\) (in the conventional SK sign conventions).
Integrating \(t\) from \(0\) to \(1\) compares the endpoint free energies and yields
\(\limsup_{N\to\infty}\FN(\beta,h)\le \mathcal{P}_k(m,q)\),
and then optimizing over \((m,q)\) gives \(\limsup \FN \le \mathcal{P}(\xi,h)\).

\paragraph{Why Guerra's discrete functional matches Parisi's PDE formulation.}
Guerra's \(\mathcal{P}_k(m,q)\) is defined by a finite recursion: at each step one applies a log-moment transform to a Gaussian increment,
with the parameters \(m_\ell\) playing the role of ``inverse temperatures'' for these intermediate transforms.
As the grid is refined, this recursion converges to the solution of the Parisi PDE (a nonlinear Hamilton--Jacobi--Bellman equation),
which is precisely the analytic form Parisi originally used.
We return to the PDE representation in \Cref{sec:parisi-measures}, but conceptually the key point is:
\begin{quote}
Guerra's recursion is a discrete-time version of Parisi's continuous-time PDE.
\end{quote}
This is one of the cleanest bridges between the physics ansatz (written as a functional order parameter) and a rigorous inequality.

\subsection{Aizenman--Sims--Starr's extended variational principle}\label{sec:ASS}

Aizenman--Sims--Starr (AS$^2$) developed a cavity-driven variational framework \cite{AizenmanSimsStarr2003}
in terms of \emph{random overlap structures} (ROSt).
Very roughly, one considers an abstract random array of overlaps together with random ``weights''
and defines a functional that mimics the free energy increment when adding a new spin.
The AS$^2$ principle states that the thermodynamic free energy can be expressed (or bounded sharply)
as an infimum of this functional over all admissible ROSt.

This perspective clarifies why Ruelle probability cascades appear naturally:
they are canonical candidates---and, in appropriate senses, extremal objects---in a variational problem over overlap structures.

\paragraph{How AS$^2$ connects to Guerra--Talagrand.}
One can view Guerra's RSB bound as evaluating the AS$^2$ variational functional on a particular hierarchical ROSt,
essentially the one generated by a (finite-depth) Ruelle cascade consistent with the chosen grid \((m,q)\).
Talagrand's 2006 theorem can then be read as showing that no admissible ROSt can do better than the Parisi value:
the cascade ROSt is not only a good candidate, but (in the relevant sense) optimal.
Thus AS$^2$ provides an organizing ``variational envelope'' around the Guerra bound and helps explain why cascades are natural objects on both sides of the Parisi formula.

% ============================================================
\section{The Parisi functional and Talagrand's 2006 theorem}\label{sec:parisi-functional}

With the thermodynamic limit in hand (\Cref{thm:thermo-limit}) and Guerra's interpolation bound (\Cref{sec:guerra-bound}),
the Parisi formula becomes the problem of identifying the limiting free energy with the same variational quantity.
Guerra's argument supplies the upper bound
\(
\limsup_{N\to\infty}\FN(\beta,h)\le \mathcal{P}(\xi,h).
\)
Talagrand's 2006 breakthrough is the matching lower bound, thereby establishing the Parisi variational principle as an equality for the SK model and for a broad class of mixed \(p\)-spin models.

We have already defined the discrete Parisi functional and stated Guerra's upper bound in \Cref{sec:guerra-bound}.
Talagrand's contribution is not a formal reversal of Guerra's interpolation.
Rather, the proof extracts sufficient quantitative control and structural information about the Gibbs measure to force the Parisi recursion to emerge as the correct effective description of the free energy.

\begin{theorem}[Parisi formula for SK and an even/convex class \cite{Talagrand2006Parisi}]\label{thm:talagrand-parisi}
For the SK model (and, more generally, mean-field Gaussian spin glasses with covariance \(N\xi(R_{1,2})\)
under suitable convexity/regularity assumptions on \(\xi\)), the limiting free energy exists and equals the Parisi variational value:
\[
\lim_{N\to\infty}\frac1N\,\E\log \sum_{\sigma\in\Sigma_N}
\exp\Bigl(-\HN(\sigma)+h\sum_{i=1}^N\sigma_i\Bigr)
=
\mathcal{P}(\xi,h).
\]
\end{theorem}

\paragraph{A schematic proof sketch.}
Guerra's interpolation yields \(\limsup \FN \le \mathcal{P}(\xi,h)\).
Talagrand proves the complementary inequality \(\liminf \FN \ge \mathcal{P}(\xi,h)\) by showing that the variational bound is accompanied by a forced structural mechanism.

Fix a discrete Parisi grid \((m,q)\).
Talagrand considers coupled systems of multiple replicas whose mutual overlaps are constrained to lie in narrow windows around the prescribed values \(q_\ell\).
He then studies the corresponding constrained free energy through two intertwined steps:
\begin{itemize}
\item \emph{A cavity step} comparing the \(N\)-spin and \((N-1)\)-spin systems, arranged so that the free energy increment resembles the recursion that defines the discrete Parisi functional.
\item \emph{An interpolation step} (with Gaussian integration by parts) used to control the effect of the overlap constraints and to bound the error terms uniformly along the recursion.
\end{itemize}
A central difficulty is to show that, after optimizing over constraints, the restriction to the prescribed overlap windows does not reduce the free energy at leading order.
This is precisely where the quantitative inequality technology developed in Talagrand's earlier work becomes essential: it supplies estimates strong enough to keep the constrained system aligned with the Parisi recursion across multiple levels.
Finally, one optimizes over \((m,q)\) and refines the grid, yielding \(\liminf \FN \ge \mathcal{P}(\xi,h)\) and hence the equality in \Cref{thm:talagrand-parisi}.

\paragraph{After 2006: extensions of the Parisi formula.}
Talagrand's theorem covers SK and a broad even/convex class.
Panchenko later extended the Parisi formula to general mixed \(p\)-spin models, including odd \(p\), by a robust argument that preserves the same variational principle across the natural model class \cite{PanchenkoParisi2014}.

% ============================================================
\section{Parisi measures: from a formula to an optimizer}\label{sec:parisi-measures}

The Parisi formula identifies the limiting free energy with the value of a variational problem.
But the variational problem comes with an optimizer, and understanding that optimizer is where the formula begins to speak about \emph{structure}.
Talagrand repeatedly emphasized that one should not stop at the number \(\mathcal{P}(\xi,h)\): the minimizing measure is a rigorous order parameter,
and its properties (support, atoms, and dependence on \((\xi,h)\)) encode the distinction between replica symmetry and symmetry breaking, as well as finer phase information.

This section records the definition of the Parisi measure and briefly explains Talagrand's PDE-based viewpoint, which makes the minimizer accessible to analysis.

\subsection{From discrete grids to measures: \(\mathcal{P}(\xi,h;\mu)\)}\label{sec:parisi-grids-to-measures}
In \Cref{sec:parisi-functional} we defined the Parisi variational problem through the discrete Guerra--Parisi functional
\((m,q)\mapsto \mathcal{P}_k(m,q)\), and wrote
\[
\mathcal{P}(\xi,h)=\inf_{k,(m,q)} \mathcal{P}_k(m,q).
\]
There is an equivalent continuum formulation in which one minimizes over probability measures \(\mu\) on \([0,1]\).
To connect the two, note that any discrete grid \((m,q)\) with
\(
0=q_0<q_1<\cdots<q_k<q_{k+1}=1
\)
and
\(
0=m_0\le m_1\le \cdots\le m_k\le m_{k+1}=1
\)
defines a nondecreasing right-continuous step function \(\alpha:[0,1]\to[0,1]\) by
\[
\alpha(q)=m_\ell \quad \text{for } q\in [q_\ell,q_{\ell+1}),
\]
which is the cumulative distribution function of a probability measure \(\mu\) on \([0,1]\) via \(\alpha(q)=\mu([0,q])\). One can write $\mathcal{P}_k(m,q)$ as $\mathcal{P}(\xi, h;\mu)$. It turns out that the map $\mu \mapsto \mathcal{P}(\xi,h;\mu)$ extends by continuity to all probability measures supported on $[0,1]$. Given \(\mu\), one can approximate its CDF \(\alpha\) by step functions (equivalently, approximate \(\mu\) by atomic measures),
so minimizing over \((m,q)\) at finer and finer grids is the same as minimizing over \(\mu\).

\begin{definition}[Parisi measure]\label{def:parisi-measure}
A \emph{Parisi measure} (for given \(\xi,h\)) is a probability measure \(\mu\) on \([0,1]\) that minimizes the Parisi functional
\(\mu\mapsto \mathcal{P}(\xi,h;\mu)\).
(Equivalently, one may minimize over discrete grids \((m,q)\), which encode \(\mu\) through the associated step-function CDF \(\alpha\); see above.)
When the minimizer is unique (see \Cref{thm:unique-parisi}), we denote it by \(\mu_{\mathrm{P}}\) and call it \emph{the} Parisi measure.
\end{definition}

\subsection{Talagrand's work on Parisi measures: the PDE viewpoint and Parisi's original formulation}

In \cite{Talagrand2006ParisiMeasures}, Talagrand studied analytic properties of the Parisi functional and its minimizers
(differentiability, variational characterizations, and the interpretation of the minimizer as an ``order parameter'' in a rigorous sense).
A central representation of the Parisi functional is via a nonlinear PDE (equivalently, a stochastic control problem), which is close in spirit to Parisi's original physics formulation.

To describe it, let \(\mu\) be a probability measure on \([0,1]\) and set \(\alpha(q)=\mu([0,q])\).
Define \(\Phi(q,x)\) as the solution of the Parisi PDE
\begin{equation}\label{eq:ParisiPDE}
\partial_q \Phi(q,x)
=
-\frac{\xi''(q)}{2}\Big(\partial_{xx}\Phi(q,x)+\alpha(q)\big(\partial_x\Phi(q,x)\big)^2\Big),
\qquad
\Phi(1,x)=\log\cosh(x).
\end{equation}
Then the Parisi functional at \(\mu\) is
\begin{equation}\label{eq:ParisiFunctionalMu}
\mathcal{P}(\xi,h;\mu)
=
\log 2+\Phi(0,h)
-\frac12\int_0^1 q\,\xi''(q)\,\alpha(q)\,\dd q,
\end{equation}
and the Parisi value is \(\mathcal{P}(\xi,h)=\inf_{\mu}\mathcal{P}(\xi,h;\mu)\). These are the same quantities that we defined previously using a different mechanism. 
In particular, the discrete Guerra functional \(\mathcal{P}_k(m,q)\) from \Cref{def:discrete-parisi} is recovered when \(\alpha\) is piecewise constant,
and~\eqref{eq:ParisiPDE} is viewed as the continuum limit of the discrete recursion for \(X_\ell\).
This is the sense in which Parisi's original ``functional order parameter'' formulation and Guerra's discrete interpolation functional are the same object, written at different resolutions.

Talagrand's contribution in \cite{Talagrand2006ParisiMeasures} was to treat the minimizer \(\mu_{\mathrm P}\) as a concrete mathematical object:
the PDE representation yields variational/stationarity conditions that determine where \(\mu_{\mathrm P}\) places mass,
and it allows one to study how the optimizer depends on \((\xi,h)\), thereby connecting the free energy formula to phase structure questions.

\subsection{Later developments: uniqueness}

Once the Parisi formula is established, it becomes meaningful to ask whether the minimizer is \emph{unique},
how it depends on parameters, and what its support says about phases (RS vs RSB).
Auffinger and Chen proved uniqueness of the Parisi minimizer \cite{AuffingerChenUnique2015},
which is a cornerstone for this analytic ``second wave'' of rigor:
it ensures that the order parameter is well-defined, stable under perturbations,
and can be studied with PDE/variational tools in a way that truly reflects the model (rather than an artifact of non-uniqueness).

\begin{theorem}[Uniqueness of the Parisi minimizer \cite{AuffingerChenUnique2015}]\label{thm:unique-parisi}
For admissible \(\xi\) and external field \(h\), the Parisi functional has a unique minimizer (a unique Parisi measure).
\end{theorem}

\subsection{The AT line and the RS region}\label{sec:AT-line}

A persistent analytic challenge is to identify the replica-symmetric region sharply.
For SK, a typical RS stability condition involves the solution \(q\) of \eqref{eq:RS-fixed-point} and the AT inequality
\begin{equation}\label{eq:AT-condition}
\beta^2\,\E\sech^4(\beta\sqrt{q}\,Z+h)\le 1,
\qquad Z\sim\mathcal{N}(0,1).
\end{equation}
Talagrand emphasized (already in 2006, and systematically in his books) that identifying the exact RS region is subtle:
it requires extracting sharp information from the Parisi minimizer, not merely proving the Parisi formula.

Rigorous results that connect to the AT line include:
\begin{itemize}
\item Guerra's 2006 analysis via quadratic replica coupling and Parisi variational methods \cite{GuerraAT2006};
\item Chen's 2021 verification of RS above the AT line in the SK model with centered Gaussian external field \cite{ChenAT2021};
\item Brennecke--Yau's 2022 enlarged RS region obtained by a refined Morita/second-moment strategy \cite{BrenneckeYau2022};
\item Mourrat--Schertzer's 2026 mixed \(p\)-spin counterexamples, showing that AT-type predictions need not match the true RS/RSB transition in general \cite{MourratSchertzer2026}.
\end{itemize}
These results show that AT-type criteria are not universally exact across all mixed \(p\)-spin models.
For the classical SK model, whether the AT condition \eqref{eq:AT-condition} gives the exact RS/RSB boundary remains a subtle frontier,
and Talagrand repeatedly emphasized identifying this sharp boundary as a central open problem.\footnote{In \cite{Talagrand2007Obnoxious}, Talagrand refers to a collection of such questions as “obnoxious problems.”}

% ============================================================
\section{Spherical models and parallel variational principles}\label{sec:spherical}

Alongside the Ising SK model, spherical mean field spin glasses provide a second ``laboratory'' in which Parisi-type variational principles can be proved and analyzed.
The advantage is largely analytic: the configuration space is a smooth manifold, so several steps (e.g., interpolation, convexity, and regularity arguments) become cleaner,
while the key disorder-driven objects---Gibbs measures, overlaps, and hierarchical organization---remain central.

\paragraph{The spherical mixed \(p\)-spin model.}
In the spherical setting, spins live on the Euclidean sphere
\[
\mathbb{S}^{N-1}(\sqrt{N})=\{\sigma\in\R^N:\norm{\sigma}_2=\sqrt{N}\},
\]
and the Hamiltonian is a centered Gaussian field with covariance
\[
\E\HN(\sigma^1)\HN(\sigma^2)=N\,\xi\Big(\frac{\inner{\sigma^1}{\sigma^2}}{N}\Big),
\]
so the overlap is the normalized inner product \(R_{1,2}=\inner{\sigma^1}{\sigma^2}/N\).

\paragraph{The spherical Parisi/Crisanti--Sommers formula.}
Physically, Crisanti and Sommers derived the corresponding variational principle for the free energy in 1992 \cite{CrisantiSommers1992}.
Talagrand's 2006 paper \cite{Talagrand2006Spherical} gave a rigorous proof of the spherical Parisi/Crisanti--Sommers formula
for a broad (even) class of mixtures, using methods parallel to the Ising case but often cleaner analytically.

\paragraph{Later developments.}
Subsequent work extended and refined this picture; for example, Chen proved the Parisi formula for spherical mixed \(p\)-spin models in generality in 2013 \cite{ChenSpherical2013},
and later papers studied fine properties (regularity, structure of minimizers, fluctuations) that are comparatively easier to access in the spherical world.

% ============================================================
\section{From identities to geometry: overlap structure, ultrametricity, and cascades}\label{sec:identities-geometry}

This ``geometry'' section is one of the places where Talagrand's post-Parisi program is most visible.
After the free energy is identified, the next step is to understand \emph{what the Gibbs measure looks like}.
Talagrand's 2008--2010 pure-state construction (see \Cref{sec:pure-states}) uses extended Ghirlanda--Guerra identities as a key hypothesis,
and Panchenko's ultrametricity theorem turns these same identities into a global structural statement about overlap arrays.
In that sense, overlap identities provide the bridge from thermodynamics (free energy) to geometry (states and hierarchy).

\paragraph{Ultrametricity: Panchenko's theorem.}
A major milestone is Panchenko's proof of Parisi ultrametricity under overlap identities \cite{PanchenkoUltrametricity2013}.
Informally: if the limiting overlap array satisfies the Ghirlanda--Guerra identities (which one can enforce via perturbations),
then the overlap structure is ultrametric.

\begin{theorem}[Ultrametricity from Ghirlanda--Guerra (informal) \cite{PanchenkoUltrametricity2013}]
Under appropriate hypotheses ensuring the Ghirlanda--Guerra identities, the limiting overlap array is almost surely ultrametric:
for three replicas, \(R_{1,2}\ge \min\{R_{1,3},R_{2,3}\}\) asymptotically.
\end{theorem}

\paragraph{Cascades as canonical ultrametric measures.}
Ultrametric overlap structures are naturally encoded by Ruelle probability cascades \cite{Ruelle1987,BolthausenSznitman1998}.
These are random hierarchical measures on an infinite tree whose leaf weights have Poisson--Dirichlet statistics,
and whose overlap of two leaves is determined by the depth of their most recent common ancestor.

The connection to earlier sections is twofold.
First, cascades are the limiting Gibbs measures of solvable hierarchical models (REM/GREM), so they provide explicit examples where the physics picture is exact (\Cref{sec:rem-grem}).
Second, cascades are built into Guerra's interpolation: the hierarchical reference system used in the RSB bound is essentially a finite-depth cascade consistent with the grid \((m,q)\) (\Cref{sec:guerra-bound}).
Thus, once overlap identities and the Parisi variational principle are in place, the emergence of cascade structures is conceptually natural:
they are the canonical objects simultaneously compatible with (i) the thermodynamic variational problem and (ii) the hierarchical self-consistency enforced by the identities.

% ============================================================
\section{Talagrand's 2008--2010 program: constructing pure states}\label{sec:pure-states}

This section sketches Talagrand's approach to constructing pure states from overlap information, explaining the role of extended Ghirlanda--Guerra identities and an atom at the maximal overlap, and how these hypotheses led to an explicit state decomposition.

After the Parisi formula, a natural next problem was to connect overlap distributions to an explicit decomposition of the Gibbs measure.
Physics predicts that, at low temperature, the Gibbs measure splits into (infinitely many) pure states with random weights,
and that the overlap distribution records whether two replicas land in the same state or in different states.
This picture is developed extensively in \cite{MPV1987}, and it underlies the hierarchical ansatz that leads to the Parisi formula.

Talagrand made major progress on this problem, under hypotheses that are natural in the mean field setting:
extended Ghirlanda--Guerra identities and the presence of an atom at the maximal overlap.

\subsection{The heuristic: an atom at \(\qstar\) means ``same state''}

Let \(\pi\) denote the limiting distribution of \(R_{1,2}\).
If \(\pi\) has a maximal point \(\qstar\) and \(\pi(\{\qstar\})=a>0\), physics suggests that \(\qstar\)
is the typical within-state overlap and that
\[
a \approx \sum_{k\ge 1} w_k^2
\]
if \(G_N\) decomposes as a mixture of pure states with weights \(w_k\).
Talagrand's work (described below) turns this heuristic into rigorous statements in appropriate asymptotic senses.

\subsection{From a clustering mechanism to Poisson--Dirichlet weights}

The next step is to make this ``same state'' interpretation operational: one needs a mechanism that turns the near-\(\qstar\) mass of the overlap distribution into an actual decomposition of the Gibbs measure into clusters, and then to identify the law of the resulting cluster weights.

In \cite{Talagrand2008PureStatesNote}, Talagrand introduced a mechanism to extract ``clusters'' from a measure on a Hilbert ball
when there is nontrivial mass near a maximal inner product.

\begin{theorem}[Clustering lemma behind pure states (informal; cf.\ \cite{Talagrand2008PureStatesNote})]
Let \(\mu\) be a probability measure on the unit ball of a Hilbert space.
If most inner products under \(\mu^{\otimes 2}\) are bounded above by \(\qstar\) and there is nontrivial mass near \(\qstar\),
then \(\mu\) admits a decomposition into disjoint sets \(A_1,A_2,\dots\) such that within each \(A_k\),
typical inner products are close to \(\qstar\), and the near-\(\qstar\) mass of pairs is mostly explained by within-\(A_k\) pairs.
\end{theorem}

In \cite{Talagrand2010PureStates}, Talagrand combined this clustering mechanism with extended Ghirlanda--Guerra identities
to construct pure states and identify the limiting statistics of their weights.

Recall from \Cref{sec:early-anchors} that the (extended) Ghirlanda--Guerra identities~\eqref{eq:GG-schematic} are stability relations for overlaps under the Gibbs measure, and they can be enforced by a small generic perturbation of the Hamiltonian.

\begin{theorem}[Pure states from EGG + an atom at \(\qstar\) (survey form; cf.\ Theorem 2.4 in \cite{Talagrand2010PureStates})]\label{thm:purestates}
Assume the extended Ghirlanda--Guerra identities hold and the limiting overlap distribution \(\pi\) has an atom
\(\pi(\{\qstar\})=a>0\) at its maximal point \(\qstar\).
Then there exist disjoint random sets \(A_1,A_2,\dots\subset \Sigma_N\) such that:
\begin{enumerate}
\item The weight sequence \(\big(\GN(A_k)\big)_{k\ge 1}\) converges in distribution to a Poisson--Dirichlet law \(\mathrm{PD}(1-a)\).
\item Within each \(A_k\), overlaps concentrate near \(\qstar\): for every \(\varepsilon>0\), for large \(N\), with high probability,
\[
\int_{A_k\times A_k} |R_{1,2}-\qstar|\,\dd \GN(\sigma^1)\dd \GN(\sigma^2)
\le \varepsilon\,\GN(A_k)^2.
\]
\end{enumerate}
\end{theorem}

\begin{remark}
Poisson--Dirichlet laws appear naturally in hierarchical models (such as REM and GREM) and in Ruelle cascades.
Talagrand's theorem can be read as a rigorous emergence of this cascade weight structure from overlap identities and a top atom.
\end{remark}

\subsection{Atomic reduction}

Talagrand also showed that for overlap questions, each pure state can effectively be replaced by a single representative.

\begin{theorem}[Atomic reduction (survey form; cf.\ Theorem 2.5 in \cite{Talagrand2010PureStates})]\label{thm:atomic-reduction} 
Under the hypotheses of \Cref{thm:purestates}, there exists a measurable embedding \(\varphi:\Sigma_N\to B_N\subset\R^N\)
such that overlaps are well approximated by normalized inner products of \(\varphi(\sigma)\),
and \(\varphi_\#\GN\) is close to a purely atomic measure \(w_0\delta_0+\sum_{k\ge 1} w_k\delta_{x_k}\)
with \((w_k)\) asymptotically Poisson--Dirichlet.
\end{theorem}

This type of reduction is a powerful bridge from ``soft'' overlap information to ``hard'' geometric decompositions.

% ============================================================
\section{Universality, TAP, and the broader ecosystem}\label{sec:ecosystem}

The story of mean field spin glasses is not only about the Parisi formula, ultrametricity, and the geometry of overlaps.
Around these core results lies a broader ecosystem of problems and techniques that both motivated and benefited from Talagrand's work.
Some of these directions ask how robust the Parisi picture is under changes of disorder; others reinterpret SK through the TAP equations and iterative algorithms; still others study the energy landscape and extremal processes of mean field Hamiltonians.

The common thread is methodological.
Many of the advances in these neighboring areas rely on the same set of ideas that Talagrand helped sharpen into a systematic toolkit:
interpolation along carefully designed paths, concentration of measure and sharp deviation bounds, and comparison principles for Gaussian and nearly Gaussian disorder.
In this section we briefly highlight a few representative strands, emphasizing how they connect back to (or borrow from) the rigorous SK theory developed in the previous sections.

\subsection{Universality of the free energy}

Carmona and Hu proved universality results showing that the SK free energy limit is robust to non-Gaussian disorder under suitable moments \cite{CarmonaHu2006}.
Chatterjee developed additional universality and approximation methods, including a connection to Stein's method \cite{ChatterjeeStein2010}.
These results complement the Parisi formula: they show that Parisi's value often governs an entire universality class, not only the Gaussian model.

Talagrand's connection here is partly methodological:
universality proofs lean heavily on the same ``smart path'' interpolation, concentration, and quantitative inequalities
that Talagrand developed and systematized, and his books helped normalize these techniques as standard tools.

\subsection{TAP equations and iterative constructions}

The TAP approach predates the rigorous Parisi formula, and it remains important because it links thermodynamics to local fields.
Talagrand's books treat TAP equations and their role in the RS/RSB boundary as part of the core narrative of SK
(see, e.g., the SK chapter in \cite{TalagrandMF1}, which includes sections on the TAP equations and the AT line).

Bolthausen gave an iterative construction of solutions to TAP equations for SK \cite{BolthausenTAP2014},
with convergence governed by stability conditions closely related to the AT criterion.
More recently, refinements of Morita/second-moment approaches connect TAP centering to rigorous RS proofs in enlarged regions;
see, e.g., \cite{BrenneckeYau2022} and references therein.
This line of work connects spin glass theory to algorithmic ideas and to high-dimensional inference (e.g., AMP/state evolution).

\subsection{Landscapes, extremes, and complexity}

A large literature studies the energy landscape of \(p\)-spin models:
the number and distribution of critical points, the structure of near-maximizers, and extremal processes.
This includes work by many authors (e.g.,  Auffinger, Ben Arous, and \v{C}ern\'y \cite{AuffingerBenArousCerny2013} on complexity)
and later developments on the geometry and extremes of spin glass Hamiltonians.
In the spherical setting, a particularly deep strand of results initiated by Subag analyzes how the Gibbs measure concentrates near bands around critical points and how near-ground-state structure organizes at low temperature; see, e.g., \cite{Subag2017,Subag2018}.
While not Talagrand's central focus, these questions use the same probabilistic culture (Gaussian comparison, interpolation ideas, concentration and variational principles) that Talagrand helped build and codify.

\subsection{Other disordered mean field models}

Talagrand's work on random $k$-SAT and the Hopfield model \cite{Talagrand2001KSAT,Talagrand2001HopfieldCritical}
is part of a broader story: tools forged in SK (concentration, cavity, interpolation) propagate to random optimization and inference problems.
Conversely, ideas from those areas have fed back into spin glass theory and strengthened the community.

A broad overview of the models treated in Talagrand's two-volume \emph{Mean Field Models for Spin Glasses}
includes (among others) SK (with several variants), perceptron-type models, Hopfield/neural network models, diluted SK and random constraint satisfaction (e.g., \(k\)-SAT),
and related mean field disordered optimization problems; see the table of contents and notes in \cite{TalagrandMF1,TalagrandMF2}.
This breadth is part of what made the spin glass toolkit a portable set of methods across probability, optimization, and inference.

% ============================================================
\section{Books as contributions: codifying a field}\label{sec:books}

Talagrand's books are not just exposition; they shaped what became ``standard technology.''
They also illustrate a feature of Talagrand's influence that is hard to capture by citing theorems alone:
he repeatedly took techniques that were scattered across papers (often with physics intuition mixed in)
and rewrote them as stable, reusable arguments with explicit constants, clear error terms, and modular lemmas.

\subsection{2003: \emph{Spin Glasses: A Challenge for Mathematicians}}

Talagrand's 2003 monograph \cite{TalagrandBook2003} is a bridge book.
It is written for mathematicians entering the subject, and it places overlaps, cavity ideas, and inequalities at center stage.
Two lasting contributions of this book (as an object, not just as a repository) are:
\begin{itemize}
\item it makes ``overlap language'' the default coordinate system for SK;
\item it treats concentration and interpolation as infrastructure, not as ad hoc tricks.
\end{itemize}
Even when later results superseded some early bounds, the book's organization strongly influenced how the field teaches itself.

\subsection{2010--2011: \emph{Mean Field Models for Spin Glasses}, Vol.\ I--II}

The two-volume treatise \cite{TalagrandMF1,TalagrandMF2} is closer to an encyclopedia, but it is still written in Talagrand's style:
a long chain of inequalities and identities, each built to be used again.
Volume I (``Basic Examples'') develops the SK model as a first course and then treats a range of other mean field disordered models
(perceptron, Hopfield, diluted models and \(K\)-SAT, and others), emphasizing how the same core techniques reappear.
Volume II (``Advanced Replica-Symmetry and Low Temperature'') integrates the Parisi formula, Ghirlanda--Guerra identities,
and the low-temperature theory into a coherent reference.

In the same spirit, later books by other authors (e.g.,  Panchenko's \cite{PanchenkoBook2013}) helped codify the modern structural viewpoint,
but Talagrand's volumes remain a central point of entry and a record of how the subject's ``methodological spine'' was assembled.

% ============================================================
\section{A guided timeline (milestones and how they connect)}\label{sec:timeline}

To close the narrative loop, we list a streamlined timeline emphasizing how the pieces fit together.

\begin{longtable}{@{}p{1.2cm}p{0.31\linewidth}p{0.53\linewidth}@{}}
\toprule
\textbf{Era} & \textbf{Milestones} & \textbf{How they connect} \\
\midrule
\endfirsthead
\toprule
\textbf{Era} & \textbf{Milestones} & \textbf{How they connect} \\
\midrule
\endhead

1970s--1980s &
SK \cite{SK1975}, TAP \cite{TAP1977}, AT \cite{AT1978}, Parisi \cite{Parisi1979}, MPV \cite{MPV1987}, REM/GREM \cite{Derrida1981}, Ruelle cascades \cite{Ruelle1987} &
Physics proposes a hierarchical order parameter and ultrametric geometry; hierarchical toy models make the picture explicit. \\

1987--1998 &
High-\(T\) rigor \cite{AizenmanLebowitzRuelle1987,CometsNeveu1995}, non-self-averaging \cite{PasturShcherbina1991},
stochastic stability \cite{AizenmanContucci1998}, Ghirlanda--Guerra identities \cite{GhirlandaGuerra1998} &
The overlap becomes a rigorous observable with constraints; stability identities emerge as structural engines. \\

1998--2003 &
Talagrand's ``challenge'' \cite{Talagrand1998SKChallenge}, exponential inequalities \cite{Talagrand2000RSBExp},
thermodynamic limit \cite{GuerraToninelli2002}, Guerra bound \cite{Guerra2003BrokenRSB},
AS$^2$ variational principle \cite{AizenmanSimsStarr2003}, Talagrand's 2003 book \cite{TalagrandBook2003} &
The field reorganizes around interpolation/cavity and overlap language; Parisi functional becomes a rigorous bound. \\

2006 &
Talagrand proves Parisi formula \cite{Talagrand2006Parisi}; Parisi measures \cite{Talagrand2006ParisiMeasures}; spherical model \cite{Talagrand2006Spherical} &
Thermodynamics is solved; attention shifts to optimizers and structure; spherical models provide parallel confirmation. \\

2008--2015 &
Pure states \cite{Talagrand2008PureStatesNote,Talagrand2010PureStates}, ultrametricity \cite{PanchenkoUltrametricity2013},
general mixtures \cite{PanchenkoParisi2014}, uniqueness \cite{AuffingerChenUnique2015}, TAP iteration \cite{BolthausenTAP2014} &
Overlap identities and Parisi theory mature into geometric theorems (pure states, ultrametricity), and analytic refinement of minimizers. \\
\bottomrule
\end{longtable}

% ============================================================
\section{Outlook: what remains open}\label{sec:outlook}
Beyond particular theorems, Talagrand’s lasting contribution is a standard of argument and exposition: the insistence that the correct objects be identified, that the relevant inequalities be made quantitative, and that the proof architecture be robust enough to be reused. In mean field spin glasses this perspective helped transform a set of remarkably precise physical predictions into a mathematical theory with canonical order parameters, structural theorems, and a shared language. The subsequent development of the field continues to reflect that standard.

Many deep questions remain (and Talagrand repeatedly emphasized them, e.g., in \cite{Talagrand2007Obnoxious}). Here are a few of the most important ones:
\begin{itemize}
\item \textbf{Chaos and sensitivity:} how overlap structures respond to perturbations in disorder or temperature (see, e.g., \cite{ChatterjeeChaosBook,ChatterjeeDisorderChaosSK}).
\item \textbf{Sharp RS boundary:} for general mixed \(p\)-spin models, recent counterexamples show that AT-type criteria need not match the RS/RSB transition \cite{MourratSchertzer2026}; for the classical SK model, proving exact sharpness of the AT criterion remains subtle (see, e.g., \cite{GuerraAT2006,ChenAT2021,BrenneckeYau2022}).
\item \textbf{Full Gibbs measure description:} ultrametricity constrains overlaps, but turning that into complete measure-level descriptions across regimes is delicate (see, e.g., \cite{PanchenkoUltrametricity2013} and the general structural framework in \cite{PanchenkoBook2013}).
\item \textbf{Beyond mean field:} transferring insights to finite-dimensional spin glasses is a longstanding frontier; see, e.g., \cite{NewmanStein2003,BinderYoung1986}.
\end{itemize}

\end{document}